# SEQUENTIAL MULTIDIMENSIONAL SPECTRAL ESTIMATION

**RAMI KANHOUCHE**[1]

**June 14, 2005**

*Abstract*-**By considering an empirical approximation, and a new class of operators that we will call walking operators, we construct, for any positive ND-toeplitz matrix, an infinite in all dimensions matrix, for which the inverse approximates the original matrix in its finite part. A recursive hierarchical algorithm is presented for sequential dimension spectral representation. A positive comparison in calculus cost, and numerical simulation, for 2D and 3D signals, is also presented.**

*Résumé*-**En considérant une approximation empirique et une nouvelle classe des opérateurs qu'on appèle les opérateurs « Walking », on arrive à construire, pour toute matrice positive ND-Toeplitz, une matrice infinie dans toutes les dimensions, et pour laquelle l'inverse; est une approximation de la matrice originale dans sa partie finie. L'algorithme hiérarchique récursif qu'on présente produit une estimation séquentielle de la puissance spectrale multidimensionnelle. On montre aussi une comparaison positive dans le coût de calcul par rapport à la méthode de Capon, et des simulations numériques pour des signaux 2D, et 3D.**

## I Introduction

Lot of work has been done in the past in the domain of spectral analysis. As a basic mathematical problem giving solutions to many physical systems, it remains not a totally resolved problem in its multidimensional aspect. While the foundations of the theory were established by [11],[12], a very series of important advances in the field were achieved, especially, the connection to the information theory [13], by the introduction of the spectral maximum of entropy principle in [1]. Also the connection to an optimized Toeplitz linear system was made in [10] and [9]. In the two dimensional case, and more generally in the multidimensional case, actually, and contrary to the correlation matching, and maximum of entropy criteria achieved in the 1D case, two algebraic methods are used. Each got its advantage and disadvantage. The first method is what we can call the quarter-plan forced autoregressive filters [7], and the second is the Capon Minimum variance method [6],[15]. While the Capon estimation provide a good stable estimation; it got the disadvantage of high calculus coast, on the other hand the quarter-plan method follows an unreal model[2], and for that suffers from high instability with order elevation. In this text we are proposing a method that is less costly than the Capon method, and far more stable than the quarter-plan autoregressive estimation. The work done and the proposed approximation depend highly on the previous insight achieved in [8] and [14].

---

[1] PhD Student at Lab. CMLA, École Normale Supérieure de Cachan, 61, avenue du Président Wilson, 94235 CACHAN Cedex, France.Phone: +33-1-40112688, mobile: 33-6-62298219, fax: +33-1-47405901.E-mail : rami.kanhouche@cmla.ens-cachan.fr, kanram@free.fr.
Web: http://www.cmla.ens-cachan.fr/Utilisateurs/kanhouch

[2] Polynomial's roots are not forcefully inside the unit circle.

## II Multidimensional Toeplitz Characters and Spectral Block Structure

For a multidimensional signal $x(t), t \in \mathbb{Z}^d$, we will define the $\gamma$ order correlation matrix, $\gamma \in (\mathbb{Z}^+)^d$ as

$$R := E\left(X(t)\left(X(t)\right)^H\right)$$

where $X$ is defined as

$$X(t) := \begin{bmatrix} x(t_0, & t_1, & \dots & t_{d-1}) \\ x(t_0+1, & t_1, & \dots & t_{d-1}) \\ \vdots \\ x(t_0+\gamma_0-1, & t_1, & \dots & t_{d-1}) \\ x(t_0, & t_1+1, & \dots & t_{d-1}) \\ \vdots \\ x(t_0, & t_1+\gamma_1-1, & \dots & t_{d-1}) \\ \vdots \\ x(t_0, & t_1, & \dots & t_{d-1}+\gamma_{d-1}-1) \end{bmatrix}.$$

As far as this paper is concerned we will be only dealing with the case where the correlation matrix admit a toeplitz nature in all dimensions, in other words, it is $d$ time toeplitz, so it is also adequate to represent it as an operator of a $d$ dimension correlation signal $c(t)$ so that

$$c(t) \in \mathbb{C}, t \in [-\gamma, \gamma], \gamma \in (\mathbb{Z}^+)^d,$$

$$c(t) \equiv \sum_{n \in \mathbb{Z}^d} x(n+t)\, x^*(n)$$

and we will write

$R = R(c(\gamma))$ which is of size $q = \prod_{i=0}^{d-1} \gamma_i$.

It is also well known, and according to Kolmogorov's relation [12] that

$$c(t) = \frac{1}{(2\pi)^d} \oint_{|\mathbb{C}^d|=1} \left| \sum_{n \in \mathbb{Z}^d} x[n] e^{jw.(n^T)} \right|^2 e^{-jw.(t^T)}.dw.$$

By considering $\Omega$, a redistribution of the dimensions nesting, such that

$$\Omega \in \left\{ (\Omega_0, \Omega_1, \dots, \Omega_{d-1}) : \Omega_k \neq \Omega_l, \forall \Omega_l, 0 \le \Omega_l \le d-1 \right\}$$

we will define

$$q_l^{\Omega} := \begin{cases} \prod_{i=0}^{l-1} \gamma_{\Omega_i} & if \quad l \neq 0 \\ 1 & if \quad l = 0 \end{cases}$$

Of course $\Omega$ is invertible in the sense that $\Omega_k = l \Rightarrow \Omega^{-1}(l) = k$.

***Definition 1.*** *We will define the walking operator as*

$$W_{\Omega\to\Omega'}(M):[M]_{i,j}\to[M]_{X^{\Omega\to\Omega'}(i),X^{\Omega\to\Omega'}(j)}$$
$$X^{\Omega\to\Omega'}(i_0q_0^{\Omega}+i_1q_1^{\Omega}+\dots i_{d-1}q_{d-1}^{\Omega})=i_0q^{\Omega'}_{\Omega'^{-1}(\Omega_1(0))}+i_1q^{\Omega'}_{\Omega'^{-1}(\Omega_1(1))}+\dots i_{d-1}q^{\Omega'}_{\Omega'^{-1}(\Omega_1(d-1))}$$
*where* $[i_0,i_1,\dots i_{d-1}]\in(\mathbb{Z}^*)^d$ $,0\le i_k\le\gamma_{\Omega(k)}$.

From a basic point of view, the walking operator corresponds to the convention used in the relation defining the vector *X*, but, as we will see later, it represents also a very important constraint on the degree of liberty of the correlation matrix, and eventually its inverse.

***Definition 2.*** *We call for any* $t\in(\mathbb{Z}^+)^d$*, t-block matrix any matrix of size* $\prod_{i=0}^{d-1}t_i$.

A d-block matrix $M$ can be indexed by the virtue of some dimensions nesting $\Omega$ according to

$$M\left[i_0q_0^{\Omega}+i_1q_1^{\Omega}+\dots i_{d-1}q_{-1}^{\Omega}\right]\left[j_0q_0^{\Omega}+j_1q_1^{\Omega}+\dots j_{d-1}q_{-1}^{\Omega}\right],\ 0\le i_l\le\gamma_{\Omega(l)},0\le j_l\le\gamma_{\Omega(l)}.$$

***Definition 3.*** *A matrix admits* $(\gamma,u,\Omega)$*-toeplitz character when it is* $\gamma$*-block matrix and realize the following*

$$M\left[i_0q_0^{\Omega}+i_1q_1^{\Omega}+\dots i_{d-1}q_{d-1}^{\Omega}\right]\left[j_0q_0^{\Omega}+j_1q_1^{\Omega}+\dots j_{d-1}q_{d-1}^{\Omega}\right]=$$
$$M\left[i_0q_0^{\Omega}+i_1q_1^{\Omega}+(i_u-j_u)q_u^{\Omega}\dots i_{d-1}q_{d-1}^{\Omega}\right]\left[j_0q_0^{\Omega}+j_1q_1^{\Omega}+0.q_u^{\Omega}+\dots j_{d-1}q_{d-1}^{\Omega}\right]\quad if\ i_u-j_u\ge 0$$
$$M\left[i_0q_0^{\Omega}+i_1q_1^{\Omega}+(0)q_u^{\Omega}\dots i_{d-1}q_{d-1}^{\Omega}\right]\left[j_0q_0^{\Omega}+j_1q_1^{\Omega}+(j_u-i_u)q_u^{\Omega}+\dots j_{d-1}q_{d-1}^{\Omega}\right]\quad if\ i_u-j_u<0$$

***Definition 4.*** *In the case that an* $\gamma$*-block matrix contains more than one toeplitz character, such that it is* $(\gamma,u_0,\Omega),(\gamma,u_1,\Omega)\dots(\gamma,u_{s-1},\Omega)$*, toeplitz character, we will call it an* $(\gamma,u,\Omega)$*-s-toeplitz character matrix, with u equal in that case to* $[u_0,u_1,\dots u_{s-1}]$.

According to **Definition 3** we can index the matrix elements using the notation

$$M\left[i_0q_0^{\Omega}+i_1q_1^{\Omega}+\dots i_{u-1}q_{u-1}^{\Omega}\right]\left[j_0q_0^{\Omega}+j_1q_1^{\Omega}+\dots i_{u-1}q_{u-1}^{\Omega}\right](k_u)\left[i_{u+1}q_{u+1}^{\Omega}+\dots i_{d-1}q_{-1}^{\Omega}\right]\left[j_{u+1}q_{u+1}^{\Omega}+\dots j_{d-1}q_{-1}^{\Omega}\right]$$

To simplify the notations we will write it as

$$M^{\Omega}\left[b_0^{u-1}\right]\left(d_u^{u+s}\right)\left[b_{u+s+1}^{d-1}\right]\quad b_h^{h+y-1}\in(\mathbb{Z}^*)^{2y},d_h^{h+y-1}\in\mathbb{Z}^y,$$

where for $s>0$ the matrix is *s*-toeplitz characters.

***Proposition 1.*** $\forall\Omega'$*, any* $(\gamma,u,\Omega)$ *toeplitz character matrix can be transformed into* $(\gamma,u',\Omega')$ *toeplitz character using the walking transformation* $W_{\Omega\to\Omega'}(M)$*, and we have* $u'=\Omega'^{-1}(\Omega_u)$.

***Theorem 1.*** *For any full rank matrix* $(\gamma,u,\Omega)$*-s-toeplitz character matrix,* $M^{\Omega}\left[b_0^{u_0-1}\right]\left(t_{u_0},t_{u_1},\dots,t_{u_{s-1}}\right)$*, with* $t_{u_i}\in[-\infty,\infty]$*, the inverse matrix* $iM$ *is equal to*

$$iM^{\Omega}\left(l_{u_0},l_{u_1},\ldots,l_{u_{s-1}}\right)=\frac{1}{(2\pi)^s}\oint_w\left[\sum_{t_u}M^{\Omega}\left(t_{u_0},t_{u_1},\ldots,t_{u_{s-1}}\right)e^{-jw(t_u)^T}\right]^{-1}e^{jw(l_u)^T}.dw\ .$$

**Proof.** The proof is very connected to signal processing theory, with the difference of the signal being a matrix –or a block- instead of a point value. In the first part the matrix $M^{\Omega}$, can be seen as a block convolution operator. In the effort to simplify the algebra representation of this operator, we will take a continuous ordering of a half of the complex region $|\mathbb{Z}|^{s+1}=1$ as $w_i$, and with the same manner, we will consider for the time space, a discrete numbering, of the half also, as $k_i$, both $k_i$ and $w_i$ are eventually both symmetric around zero. It is well known fact that in the case where $s=d$ for $M^{\Omega}\left(t_{u_0},t_{u_1},\ldots,t_{u_{s-1}}\right)$, the convolution operator can be represented as

$$\frac{dw}{(2\pi)^d}C\lambda C^H$$

where $\lambda$ is the eigenvalues diagonal matrix defined as

$$\lambda:=\begin{bmatrix}\lambda_{w_{-\infty}} & & & & & & \\ & \ddots & & & & & \\ & & \lambda_{w_{-1}} & & & & \\ & & & \lambda_0 & & & \\ & & & & \lambda_{w_1} & & \\ & & & & & \ddots & \\ & & & & & & \lambda_{w_\infty}\end{bmatrix}$$

$$C:=\begin{bmatrix}e^{-jw_{-\infty}(k_{-\infty})^T} & \ldots & e^{-jw_{-1}(k_{-\infty})^T} & 1 & e^{-jw_1(k_{-\infty})^T} & \ldots & e^{-jw_{\infty}(k_{-\infty})^T}\\ \vdots & \ldots & & 1 & \vdots & \ldots & \vdots\\ e^{-jw_{-\infty}(k_{-1})^T} & \ldots & e^{-jw_{-1}(k_{-1})^T} & 1 & e^{-jw_1(k_{-1})^T} & \ldots & e^{-jw_{\infty}(k_{-1})^T}\\ 1 & \ldots & 1 & 1 & 1 & \ldots & 1\\ e^{-jw_{-\infty}(k_1)^T} & \ldots & e^{-jw_{-1}(k_1)^T} & 1 & e^{-jw_1(k_1)^T} & \ldots & e^{-jw_{\infty}(k_1)^T}\\ \vdots & \ldots & \vdots & 1 & & \ldots & \vdots\\ e^{-jw_{-\infty}(k_{\infty})^T} & \ldots & e^{-jw_{-1}(k_{\infty})^T} & 1 & e^{-jw_1(k_{\infty})^T} & \ldots & e^{-jw_{\infty}(k_{\infty})^T}\end{bmatrix}$$

This is can be inferred directly, as we said, from Linear Time Invariant system theory. More precisely for a given signal vector $v$ , the inverse transformation can be written as

$$x=Mv=\frac{dw}{(2\pi)^d}C\lambda C^Hv\Leftrightarrow v=M^{-1}x=\frac{dw}{(2\pi)^d}C\lambda^{-1}C^Hx\ .$$

On the other hand we can directly conclude the exactness of the previous representation by the fact that

$$\left[\frac{dw}{(2\pi)^d}CC^H\right]_{i,j} = \frac{1}{(2\pi)^d}\oint_w e^{jw(i-j)}.dw = \begin{cases}1 & if \quad i=j\\ 0 & if \quad i\neq j\end{cases}.$$

From which, we got $\frac{dw}{(2\pi)^d}CC^H = \mathbf{1}$.

In the case that $s<d$, then we will note $h := q_{d-s}^{\Omega}$ as the non-toeplitz part matrix size. As a counter part to the multidimensional Fourier transform vector in the previous definition we will define the block multidimensional Fourier transform vector as

$$F^{w_i,t} := \left[\overbrace{0 \ \dots \ 0}^{t} \ e^{-jw_i(k_{-\infty})^T} \ \overbrace{0 \ \dots \ 0}^{h-t-1} \ \dots \ \overbrace{0 \ \dots \ 0}^{t} \ e^{-jw_i(k_{\infty})^T} \ \overbrace{0 \ \dots \ 0}^{h-t-1}\right]^T,$$

$0\le t\le h-1$. By considering the Block diagonal matrix $\lambda_s$, which is by definition do contain the matrix value $\left[\sum_{k_u} M^{\Omega}\left(t_{u_0},t_{u_1},\dots,t_{u_{s-1}}\right)e^{-jw_i(t_u)^T}\right]$, at the diagonal index $w_i$, it is simple to prove that

$$\lambda_s = \frac{dw}{(2\pi)^s}C_s\left(M^{\Omega}\right)C_s^{\ H}$$

where $C_s = \left[F^{w_{-\infty},0} \quad F^{w_{-\infty},1} \quad \dots \quad F^{w_{-\infty},h} \quad \dots \quad F^{w_{\infty},h-1} \quad F^{w_{\infty},h}\right]$. Next, by observing that $C_s$ is full rank, we prove that $\lambda_s$ is full rank also, and that proves the theorem.

## III Multidimensional Extension Approximation

Our approximation that we are proposing in this depends strongly on the previous work done in the treatment of the generalized reflection coefficients [8], and [14], to better explain this extension we start with the 1 dimensional case. For a positive definite toeplitz matrix $M^N(l)$, $l := -N+1,\dots,N-1$, we define the polynomial $p^{N-1}, q^{N-1}$, as the solution for the linear system $M^N(l)\,p^{N-1} = e_0$, $M^N(l)\,q^{N-1} = e_{N-1}$. These solutions in fact can be obtained recursively according to

$$p'^N = \begin{bmatrix} p'^{N-1} \\ 0 \end{bmatrix} + \sigma_N \begin{bmatrix} 0 \\ q'^{N-1} \end{bmatrix}.$$

Where $v'$ is $v$ normalized by the first value for $p^N$, and the last for $q^N$, $\sigma$ is the famous reflection coefficient. Putting the reflection coefficient to zero is equivalent to the existence of a positive definite matrix $M^{N+1}(l)$, which coincides with $M^N(l)$, so that $M^{N+1}(l) = M^N(l), l = -N+1,\dots N-1$. By repeating the procedure until infinity. We obtain that

$$\left[M^{\infty}(l)\right]^{-1} = \rho^{-1}P.P^H \qquad (1)$$

With $P = \begin{bmatrix} \infty \nwarrow & & & & & \\ & p'^{N-1} & & & & \\ & & p'^{N-1} & & & \\ & & & \ddots & & \\ & & & & p'^{N-1} & \\ & & & & & \searrow \infty \end{bmatrix}$

And $M^N(l)\, p'^{N-1} = \rho e_0$.

In the multidimensional case, the correlation matrix $M^{\Omega}(l_0, l_1, \ldots, l_{d-1})$, $l_{d-1} := -\gamma_{\Omega_{d-1}} + 1, \ldots, \gamma_{\Omega_{d-1}} - 1$ can be always considered as a block toeplitz matrix and eventually admit the same extension as the 1 dimensional case, with the difference of $p^{N-1}$ being a matrix-vector with size $q_d^{\Omega} \times q_{d-1}^{\Omega}$, the main constraint in this case is that the new generated blocks $M^{\Omega}[\ ][\ ](l_{d-1}), l_{d-1} := -\infty, \ldots \infty$, for $|l_{d-1}| > \gamma_{\Omega_{d-1}} - 1$ are not $([\gamma_0, \ldots, \gamma_{d-2}], [\Omega_0, \ldots \Omega_{d-2}], [\Omega_0, \ldots \Omega_{d-2}])$-toeplitz character matrices. The approximation that we are proposing is to neglect this constraint and to proceed according to, in a sequential manner across each dimension. The foundation for this approximation is –as we will introduce in future work- that the block polynomial that preserve the internal toeplitz structure of the new extended matrix, do in fact coincide with the finite size polynomial in its first $q_d^{\Omega}$ line values, while for lines going from $q_d^{\Omega}$ to infinity, we observe far more less magnitude values.

***Theorem 2.*** *The multidimensional extension of the correlation matrix can be approximated in a sequential way according to the following*

$$R^{-1}\left(c\left(\left[\gamma_0, \gamma_2, \ldots \gamma_{d-x-1}, \gamma_{d-1}^{\infty}, \gamma_{d-2}^{\infty}, \ldots \gamma_{d-x}^{\infty}\right]\right)\right) = R_x^{-1}\left[b_0^{d-x-1}\right](l_{d-1}, l_{d-2}, \ldots, l_{d-x})$$

$$= \frac{1}{(2\pi)^x} \oint_{w_{d-1}^{d-x}} \begin{bmatrix} \left[\sum_{u_1^{x+}} M\left[b_0^{d-x-1}\right](u_1, \ldots, u_x) e^{jw_{d-1}^{d-x}\left(u_1^{x+}\right)^T}\right] \\ \left[\sum_{u_1^{x-1}} M\left[b_0^{d-x-1}\right](u_1, \ldots, u_{x-1}, 0) e^{jw_{d-1}^{d-x-1}\left(u_1^{x-1}\right)^T}\right]^{-1} \\ \left[\sum_{u_1^{x+}} M\left[b_0^{d-x-1}\right](u_1, \ldots, u_x) e^{jw_{d-1}^{d-x}\left(u_1^{x+}\right)^T}\right]^H \end{bmatrix} e^{-jw_{d-1}^{d-x}\left(l_{d-1}^{d-x}\right)^T} dw_{d-1}^{x} \tag{2}$$

*where*

$$M\left[b_0^{d-x-1}\right](u_1, u_2 \ldots, u_x) = R_{x-1}^{-1}\left[b_0^{d-x-1}\right](u_1, u_1 \ldots, u_{x-1})\left[(u_x, 0)\right] \tag{3}$$

$$R_{x-1}^{-1}\left[b_0^{d-x-1}\right](l_{d-1}, l_{d-2} \ldots, l_{d-x+1})\left[b_{d-x}^{d-x}\right] = W^{(0,1,\ldots,x,d-1,d-2,\ldots,d-x+1)\to(0,1\ldots,x-1,d-1,d-2,\ldots d-x+1,d-x)}\left(R_{x-1}^{-1}\left[b_0^{d-x}\right](l_{d-1}, l_{d-2} \ldots, l_{d-x+1})\right) \tag{4}$$

$$u_1^x \in \left[-\gamma_{d-1}^{\infty}, \gamma_{d-1}^{\infty}\right] \times \left[-\gamma_{d-2}^{\infty}, \gamma_{d-2}^{\infty}\right] \times \ldots \left[-\gamma_{d-x}^{\infty}, \gamma_{d-x}^{\infty}\right]$$

$$u_1^{x+} \in \left[-\gamma_{d-1}^{\infty}, \gamma_{d-1}^{\infty}\right] \times \left[-\gamma_{d-2}^{\infty}, \gamma_{d-2}^{\infty}\right] \times \ldots \left[-\gamma_{d-x+1}^{\infty}, \gamma_{d-x+1}^{\infty}\right] \times \left[0, \gamma_{d-x} - 1\right]$$

**PROOF.** Starting from $x=1$, $R_0^{-1} := \left[R\left(c\left(\left[\gamma_0,\gamma_2,\dots\gamma_{d-1}\right]\right)\right)\right]^{-1}$, relation (3) forwards the block-vector $p_x$ with columns of size $q_{d-x}^{(0,1\dots,x-1,d-1,d-2,\dots d-x+1,d-x)}$, which is next applied according to the form (1). By accepting the proposed approximation, we assume that form (1) represents the inverse of ND-Toeplitz matrix $R\left(c\left(\left[\gamma_0,\gamma_2,\dots\gamma_{d-x-1,}\gamma_{d-1}^{\infty},\gamma_{d-2}^{\infty},\dots\gamma_{d-x}^{\infty}\right]\right)\right)$, with infinite dimension according to the dimension $d-x$, and - by construction- coincide with the original matrix in its finite part. From that, and by **Theorem 1**, the inverse can be written according to the form (2). By Advancing $x$, and by **Proposition 1**, relation (4) permits us to continue the extension on the next dimension. Hence our proof is completed.

## IV Sequential Multidimensional Spectral Estimation Algorithm

**Theorem 2** presents an analytic representation of the spectral estimation and that is by taking advantage of both the **Theorem 1** and the proposed approximation. From numerical point of view, we need more finite simple Algorithm. On the contrary of the 1 dimensional case, in the general case, there is no known finite time representation for the proposed estimation, on the other hand even in the 1 dimensional case the numerical representation of the estimation is obtained always as a Fourier series of the prediction polynomial i.e. $C(jw) = \dfrac{\rho}{\left|\sum_{k=0}^{\gamma-1} p_k^{\gamma} e^{jwk}\right|^2}$

And $w$ is taken to a large finite number of points on the spectral circle, that is in an effort to approximate the continuous field of values.
By examining relation (2), and replacing by (3) starting from $x=d$, and descending in the reverse order of **Theorem 2**, we find that the spectral inverse of correlation signal $c\left(\left[\gamma_{d-1}^{\infty},\gamma_{d-2}^{\infty},\dots\gamma_0^{\infty}\right]\right)$ can be obtained by a finite series of block exponential sum, with a diminishing Block size at each step. According to the previous discussion, the algorithm takes the following form:

1. *Calculate the inverse of the correlation matrix,*
   $R_0^{-1} := \left[R^{(0,1,\dots,d-1)}\left(l_0,l_1,\dots,l_{d-1}\right)\right]^{-1}$,
   *set* $x=1$
   *set*
   $G(k) := R_0^{-1}\left[i_0 q_0^{\Omega} + i_1 q_1^{\Omega} + \dots i_{d-2} q_{d-2}^{\Omega} + k.q_{d-1}^{\Omega}\right]\left[j_0 q_0^{\Omega} + j_1 q_1^{\Omega} + \dots j_{d-2} q_{d-2}^{\Omega} + 0.q_{d-1}^{\Omega}\right]$
2. *Apply relation (1)* $M\left(w_{d-1},\dots w_{d-x}\right) = \sum_{k=0}^{\gamma_x-1} G(k)\left(w_{d-1},\dots w_{d-x+1}\right) e^{jkw_{d-x}}$
3. Set $G'\left(w_{d-1},\dots w_{d-x}\right) = M\left(w_{d-1},\dots w_{d-x}\right)\left[G(0)\right]^{-1} M^H\left(w_{d-1},\dots w_{d-x}\right)$
4. *if* $(x=d)$ *{*

*stop*
}
*else* {
*set*

$$G(k)(w_{d-1},\ldots w_{d-x}) := G'\left[i_0 q_0^{\Omega} + \ldots k.q_{d-x-1}^{\Omega}\right]\left[j_0 q_0^{\Omega} + \ldots 0.q_{d-x-1}^{\Omega}\right](w_{d-1},\ldots w_{d-x})$$

*set* $x = x+1$
*go to 2*
}

## V Numerical calculation cost

In this part we will try to shed the light on the numerical efficiency of the algorithm. For that, we will measure the number of operations needed to obtain the spectral estimation over an ND spectral grid with size $C_w := \left[C_{w_{d-1}}, C_{w_{d-2}}, \ldots C_{w_0}\right]$, where $C_{w_i}$ is the number of points in the $i$ dimension. In the following discussion we will omit the S*tep 1* from our Calculus Cost, since that $R_0^{-1}$, and more precisely $G(k), x=1$ can be obtained in different ways, like in [4], or [5], the thing that is irrelevant to the estimation method.
From the algorithm's *step 2* and *Step 3* we can write the number of operations needed at each $t := d - x + 1$ as

$$\left[\frac{3}{2}(q_{t-1})^3 + (q_{t-1}q_t)\right]\prod_{f=0}^{d-t} C_{w_{d-f-1}},$$

The total sum is equal to the sum from $t = d$, and down until $t = 1$, which can be expressed as

$$\sum_{t=1}^{d}\left(\left[\frac{3}{2}(q_{t-1})^3 + (q_{t-1}q_t)\right]\prod_{f=0}^{d-t} C_{w_{d-f-1}}\right)$$

Figure 1 graphs the above relation in function of $C_{w_0} = C_{w_1} = \ldots = C_{w_{d-1}}$, with values taken as $\gamma_0 = \gamma_1 \ldots = \gamma_{d-1} = 10$, and $d = 5$. Also for the same values Figure 2 represents the number of operations needed by Capon Estimation, which is equal to $(q_d)^2 \prod_{f=0}^{d-1} C_{w_f}$. We see clearly a very good advantage using the proposed algorithm, from a calculus coast point of view.

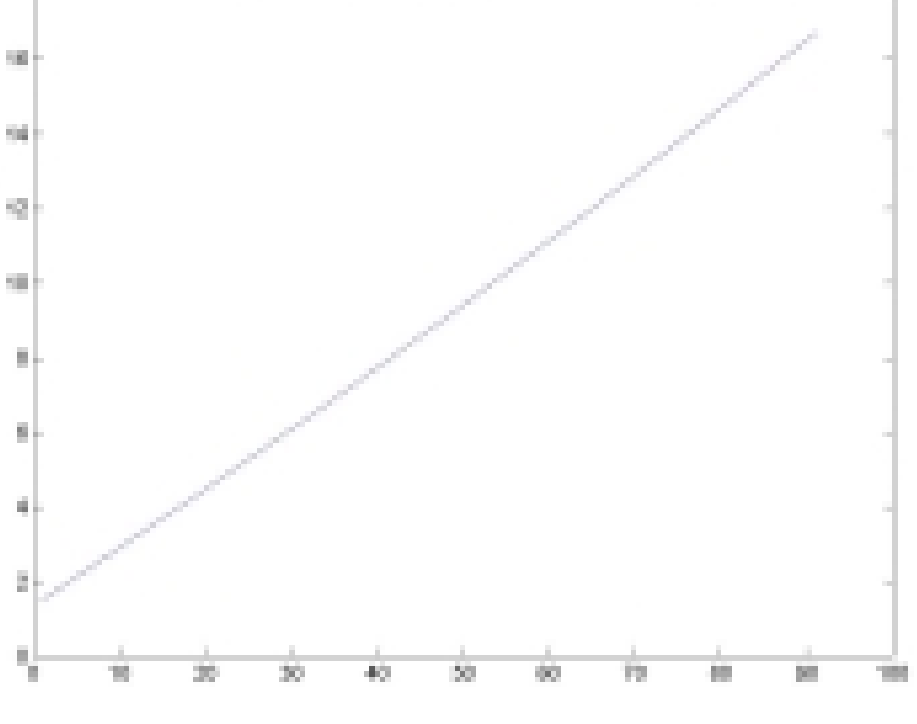

Figure 1

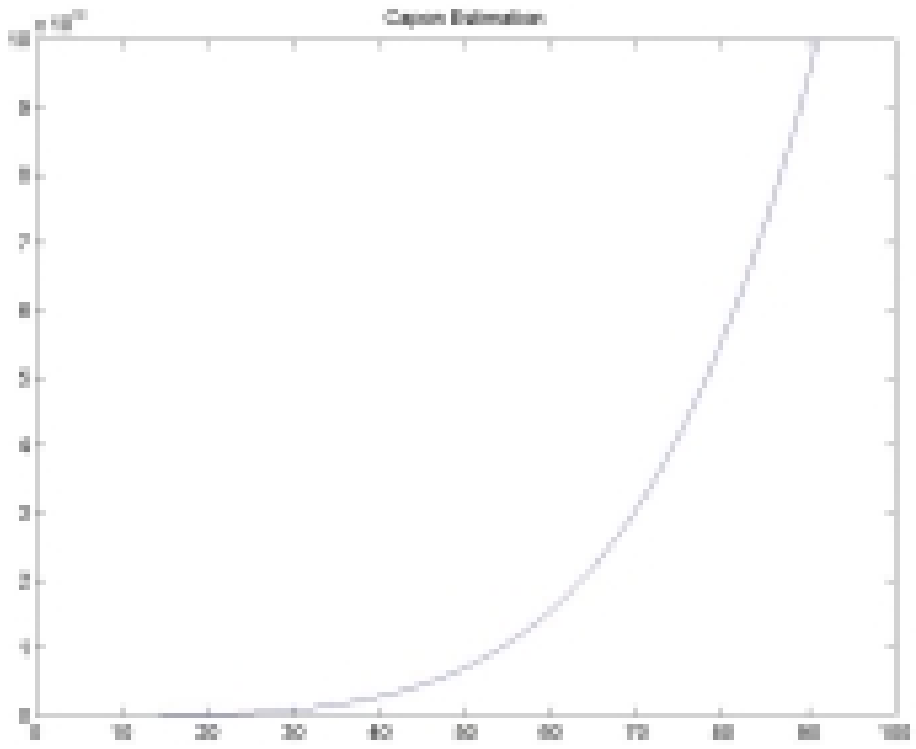

Figure 2

## VI Numerical simulations and discussion

In an effort to evaluate the proposed estimation we will present a three-dimensional test signal $c\left(w_0, w_1, w_2\right)$, defined in its spectral composition as zero everywhere except for $c\left(0.1, 0.3, 0.7\right) = 1$, $c\left(0.1, 0.6, 0.2\right) = 1$, $c\left(0.6, *, *\right) = 1$. Next, the signal was plunged in a white noise of variance 0.1. In our simulation we took our estimation on a grid of equal size in all dimensions, $C_{w_i} = 10$, and the same for the correlation matrix order, having a uniform value of $\gamma_i = 3$. Figure 3-Figure 22 show a series cuts in the cube for the simulated signal, with the noised spectrum on the left and the calculated estimation on the right.

In Figure 23, the graph representing the correlation matching accuracy of the proposed estimation for a two dimensional signal was also presented. The value present in the graph is the relative difference between the original correlation values $r_{k,l}$, and the estimated ones $\widehat{r}_{k,l}$, that is $\frac{\left|\widehat{r}_{k,l} - r_{k,l}\right|}{r_{k,l}}$.

From The presented Figures, and from other observed numerical simulations we found that the proposed algorithm provide a good estimation of the simulated signal, and contrary to the quarter-hyper plan filters proposed for the 2D case by [3][2], it enjoys a high degree of stability when augmenting in the Correlation matrix order $\gamma_i$.

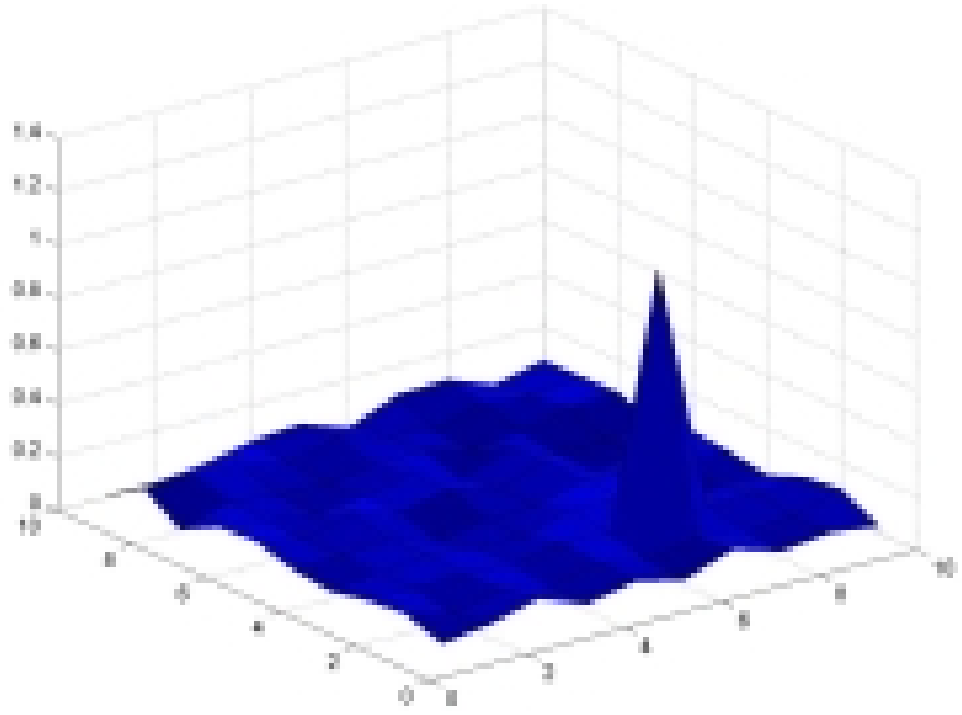

Figure 3

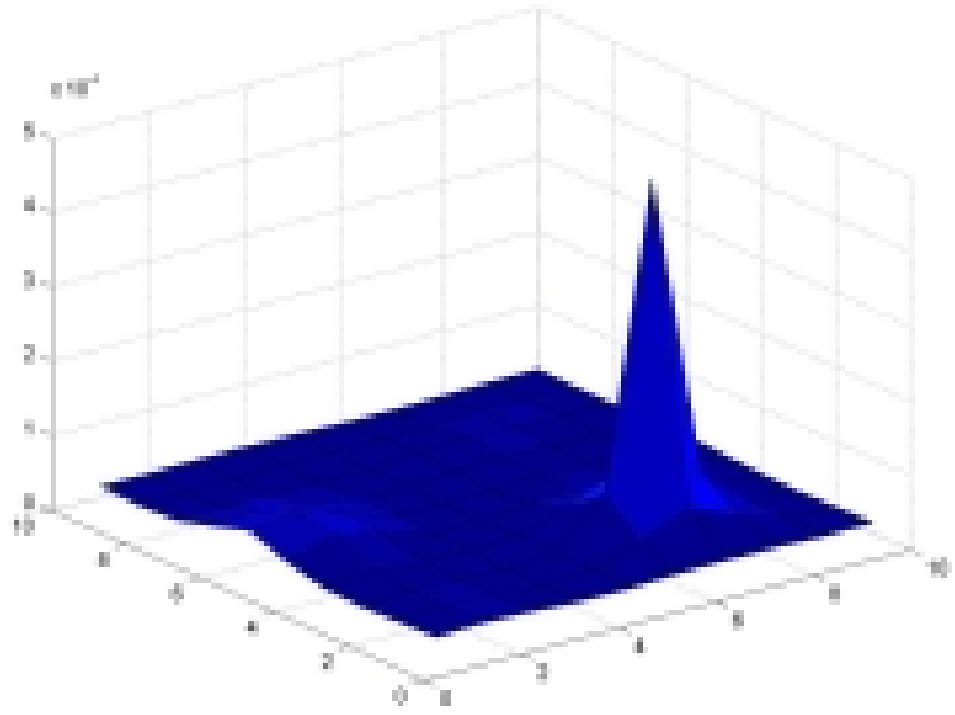

Figure 4

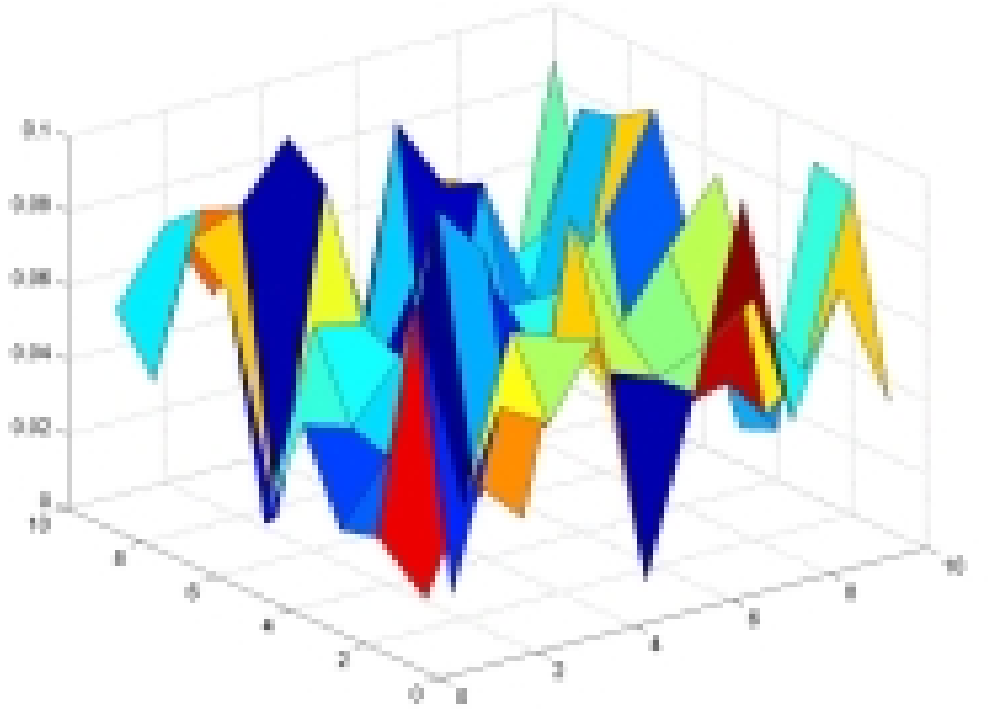

Figure 5

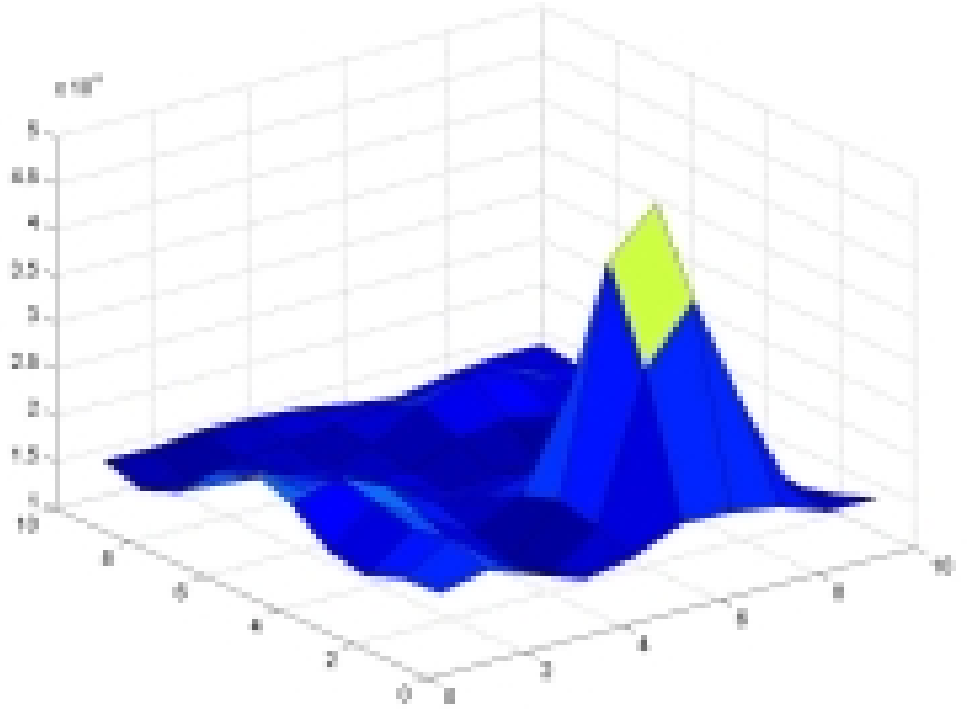

Figure 6

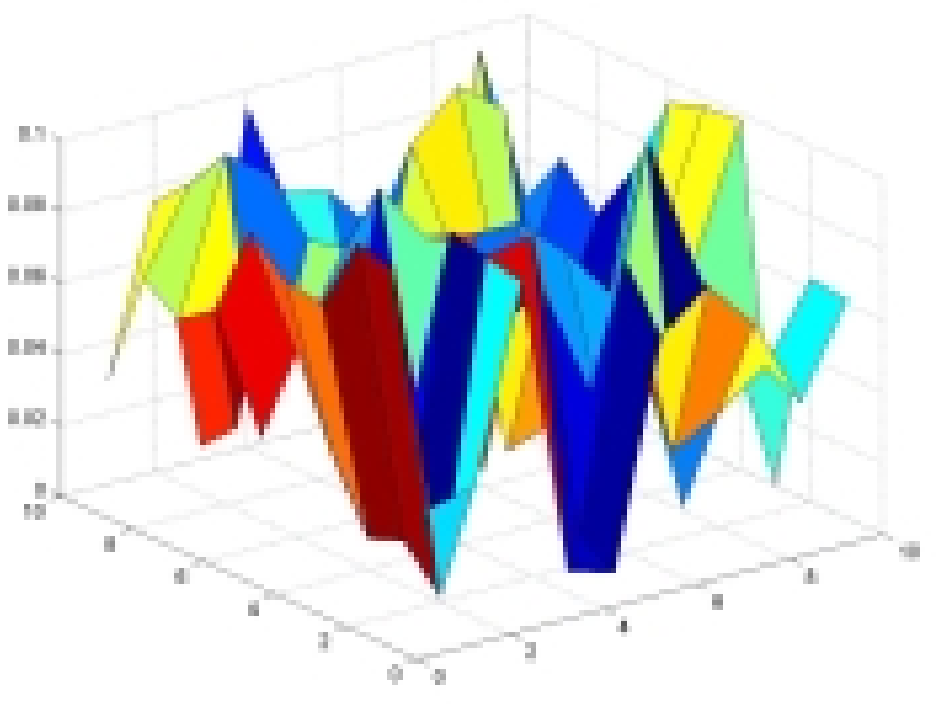

Figure 7

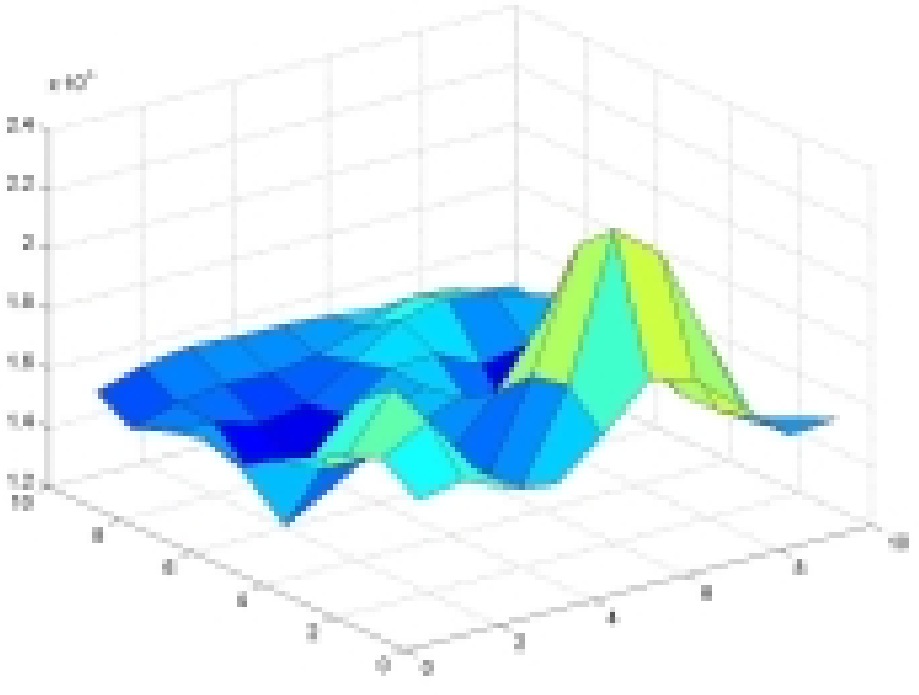

Figure 8

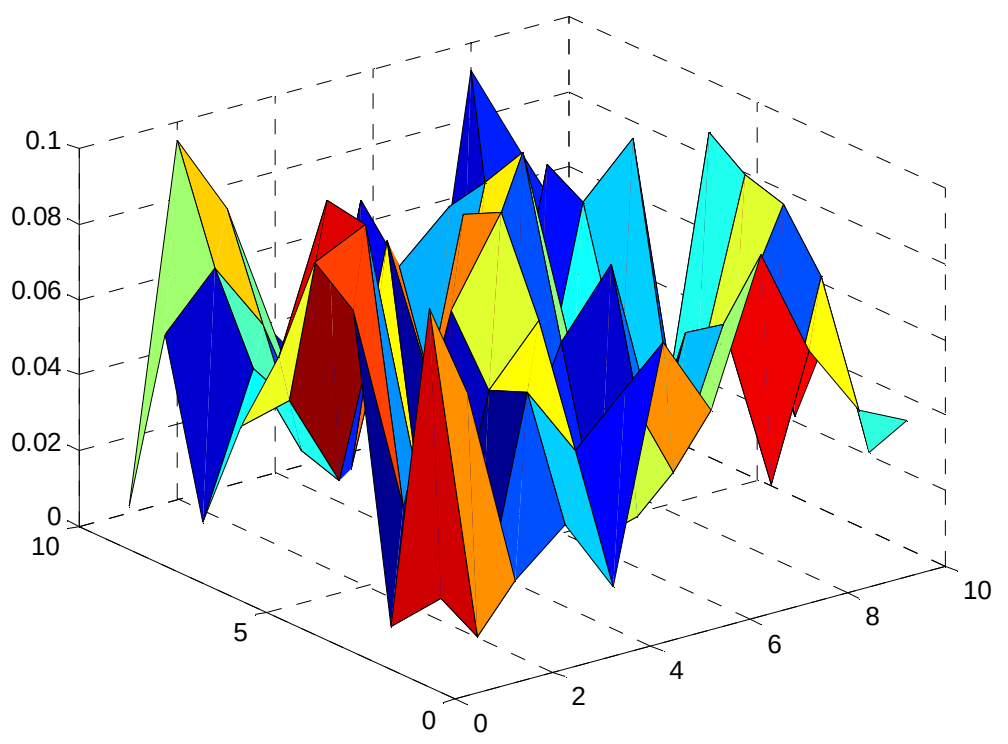

Figure 9

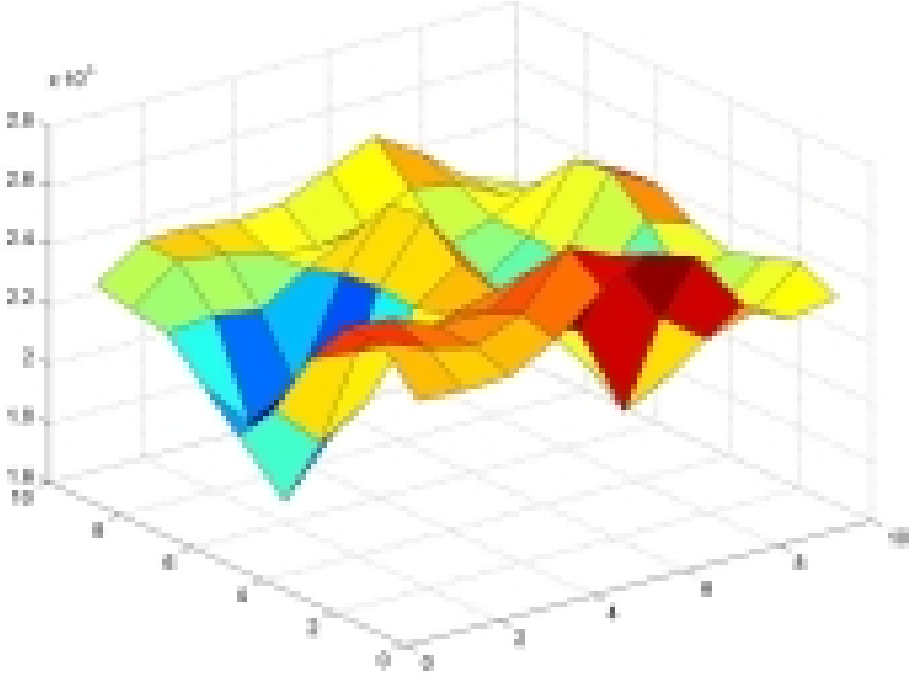

Figure 10

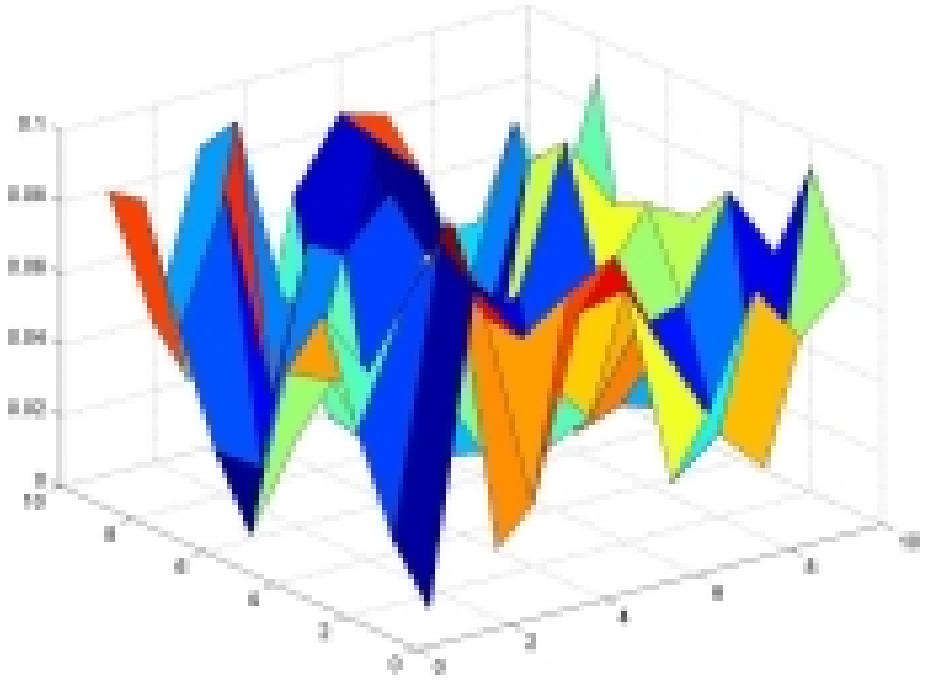

Figure 11

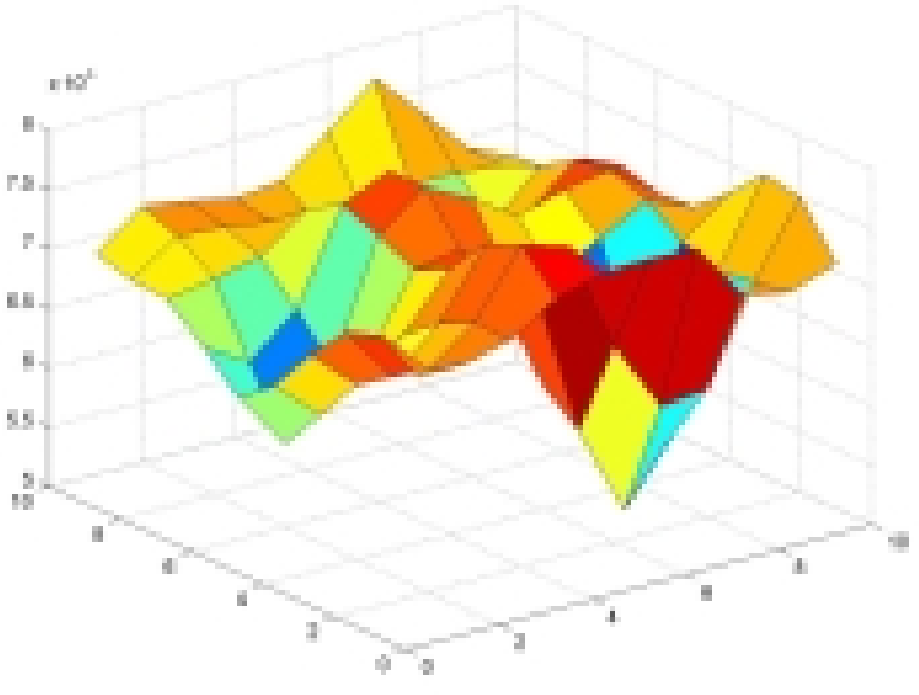

Figure 12

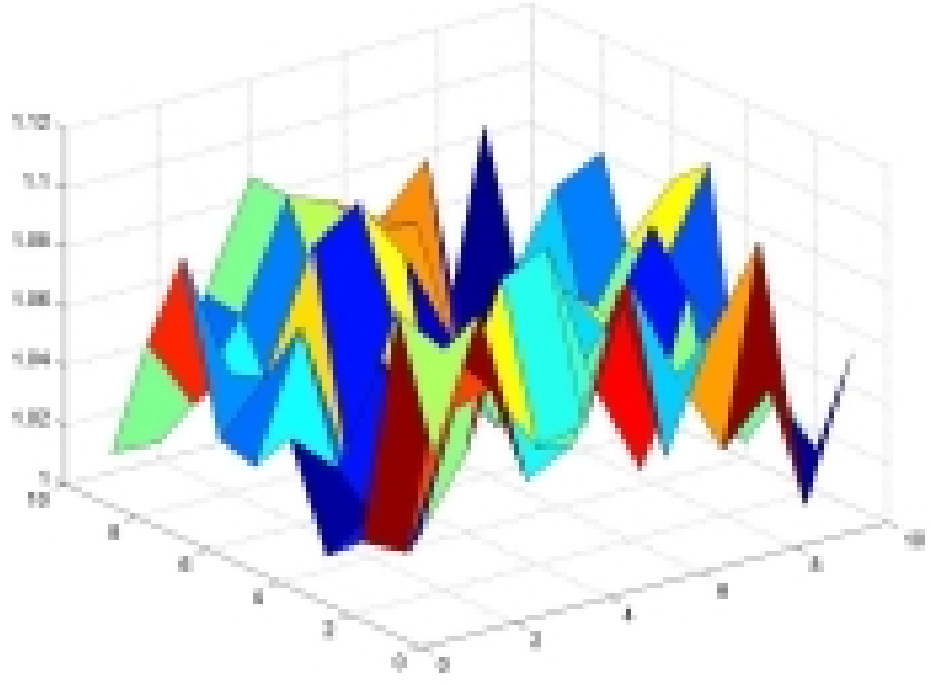

Figure 13

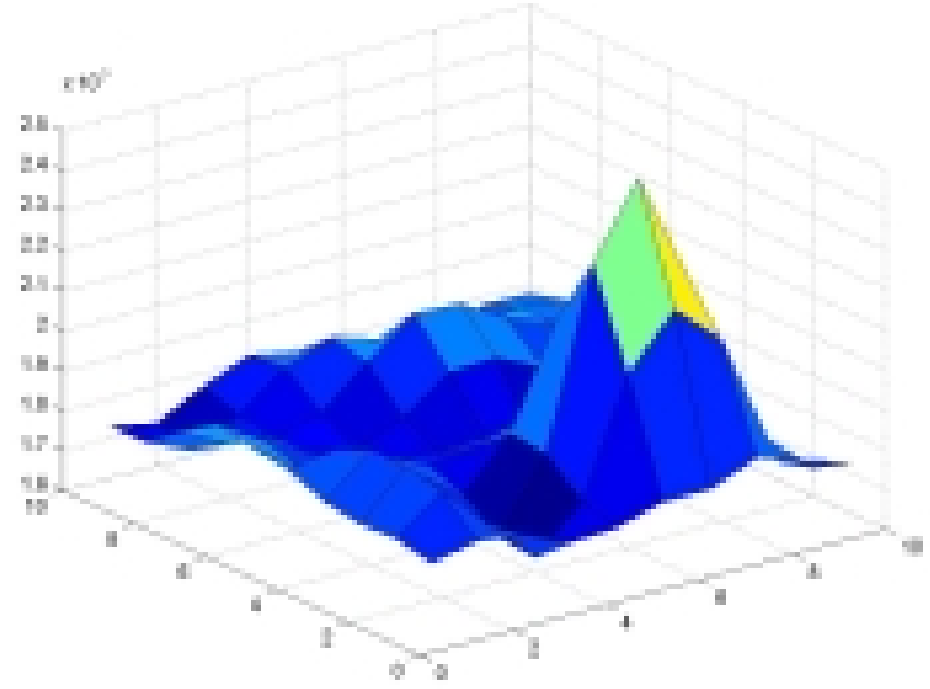

Figure 14

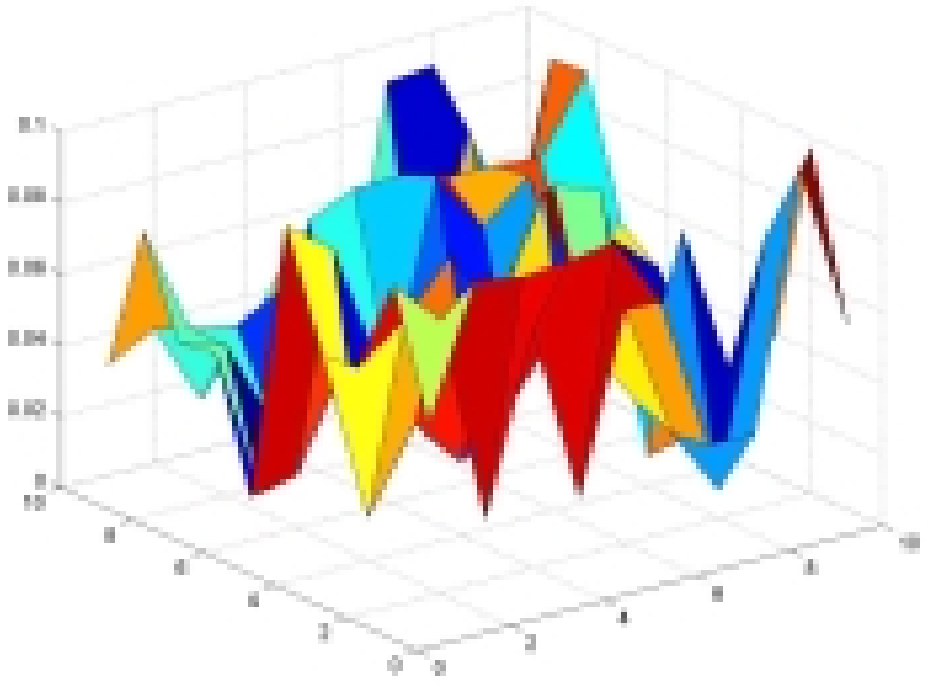

Figure 15

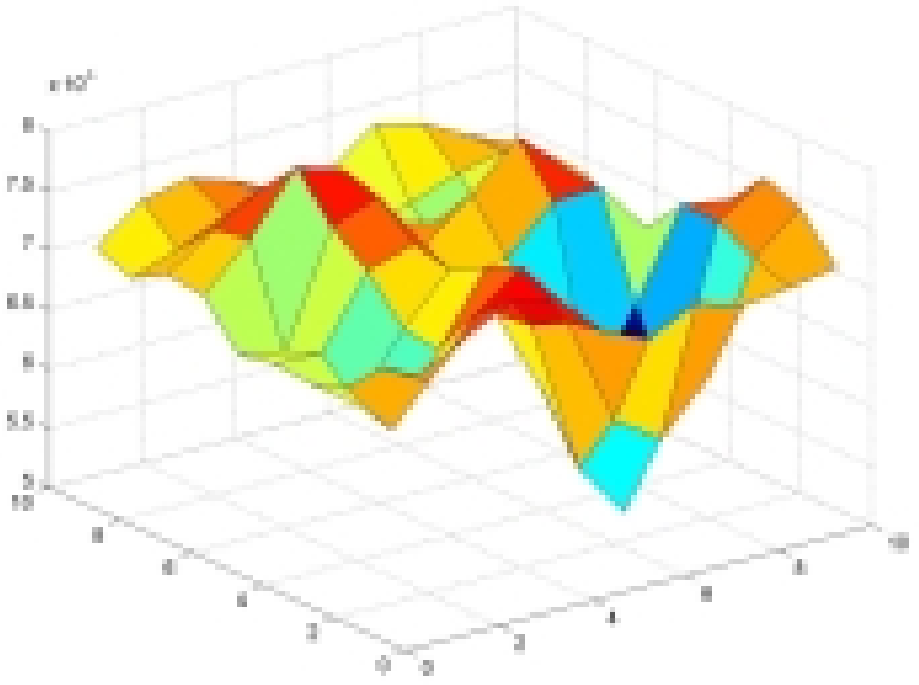

Figure 16

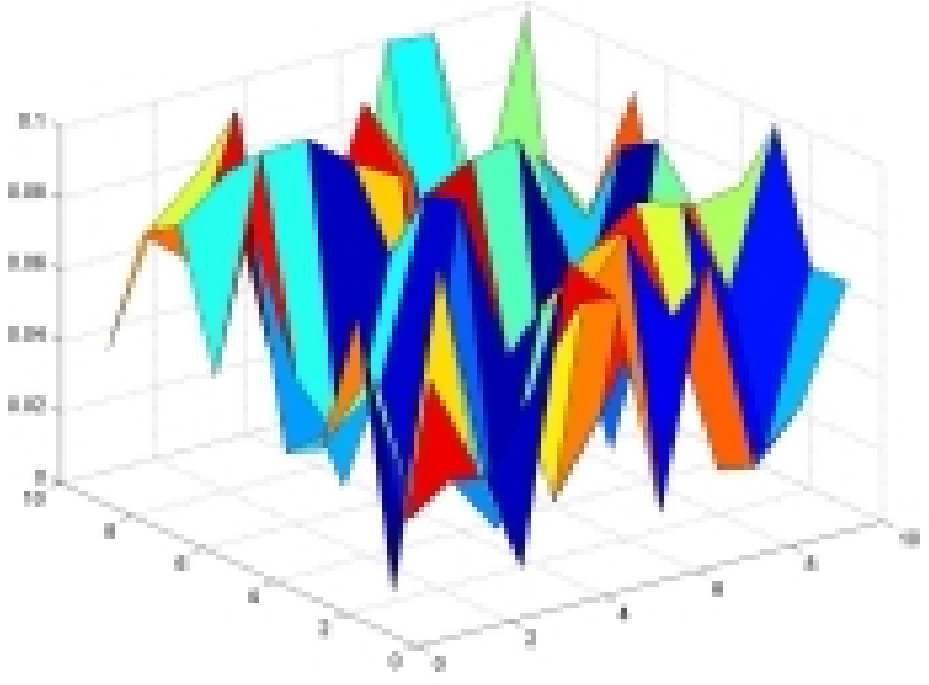

Figure 17

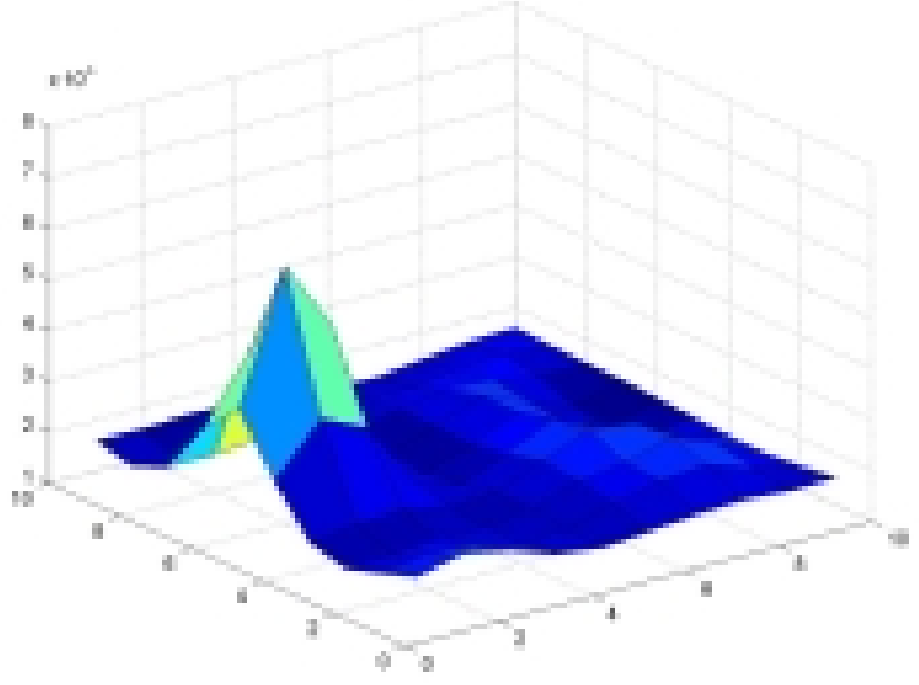

Figure 18

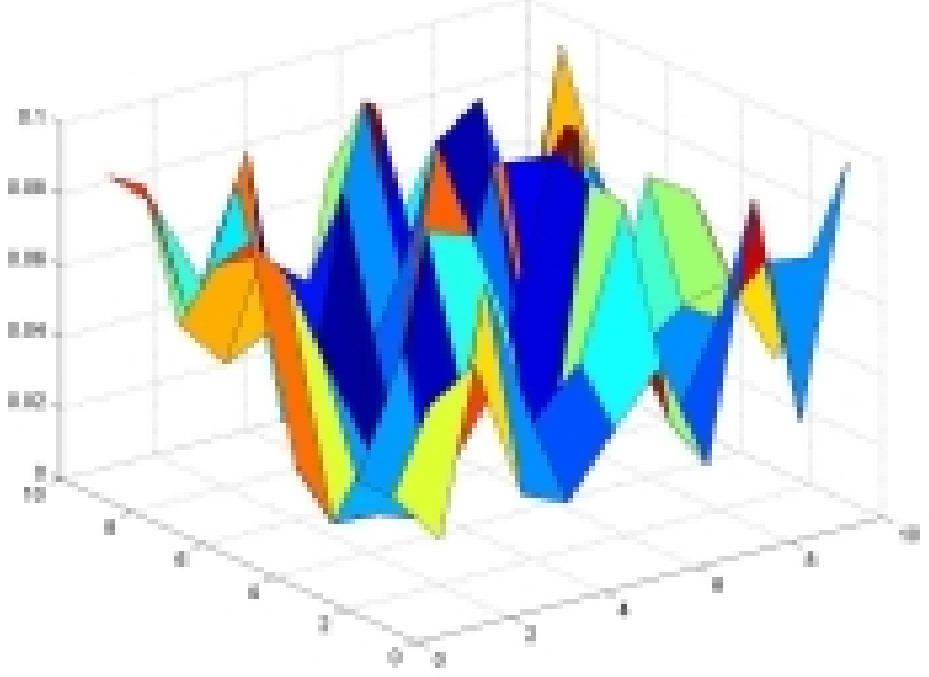

Figure 19

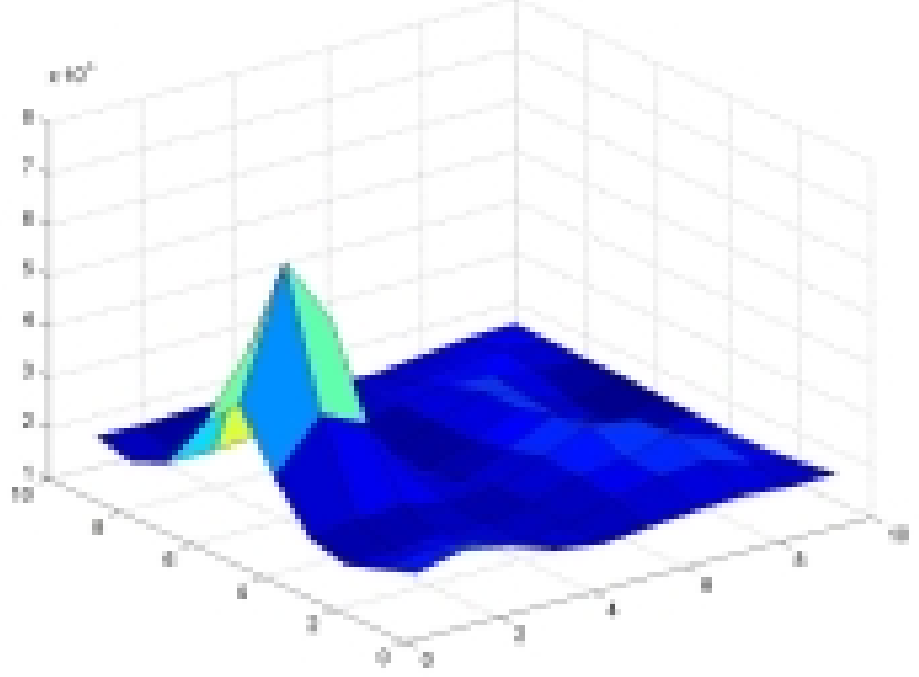

Figure 20

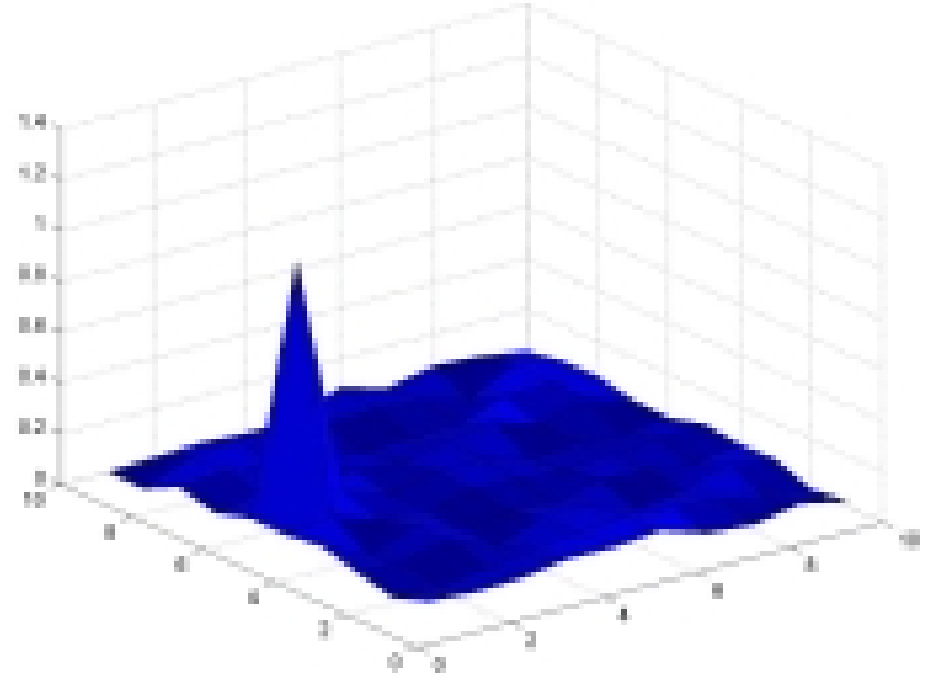

Figure 21

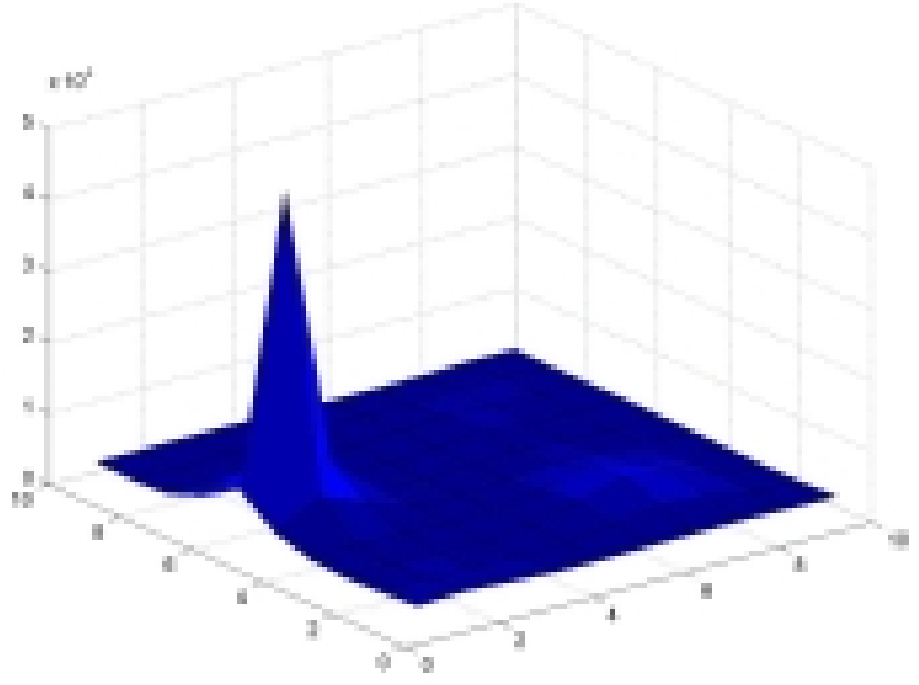

Figure 22

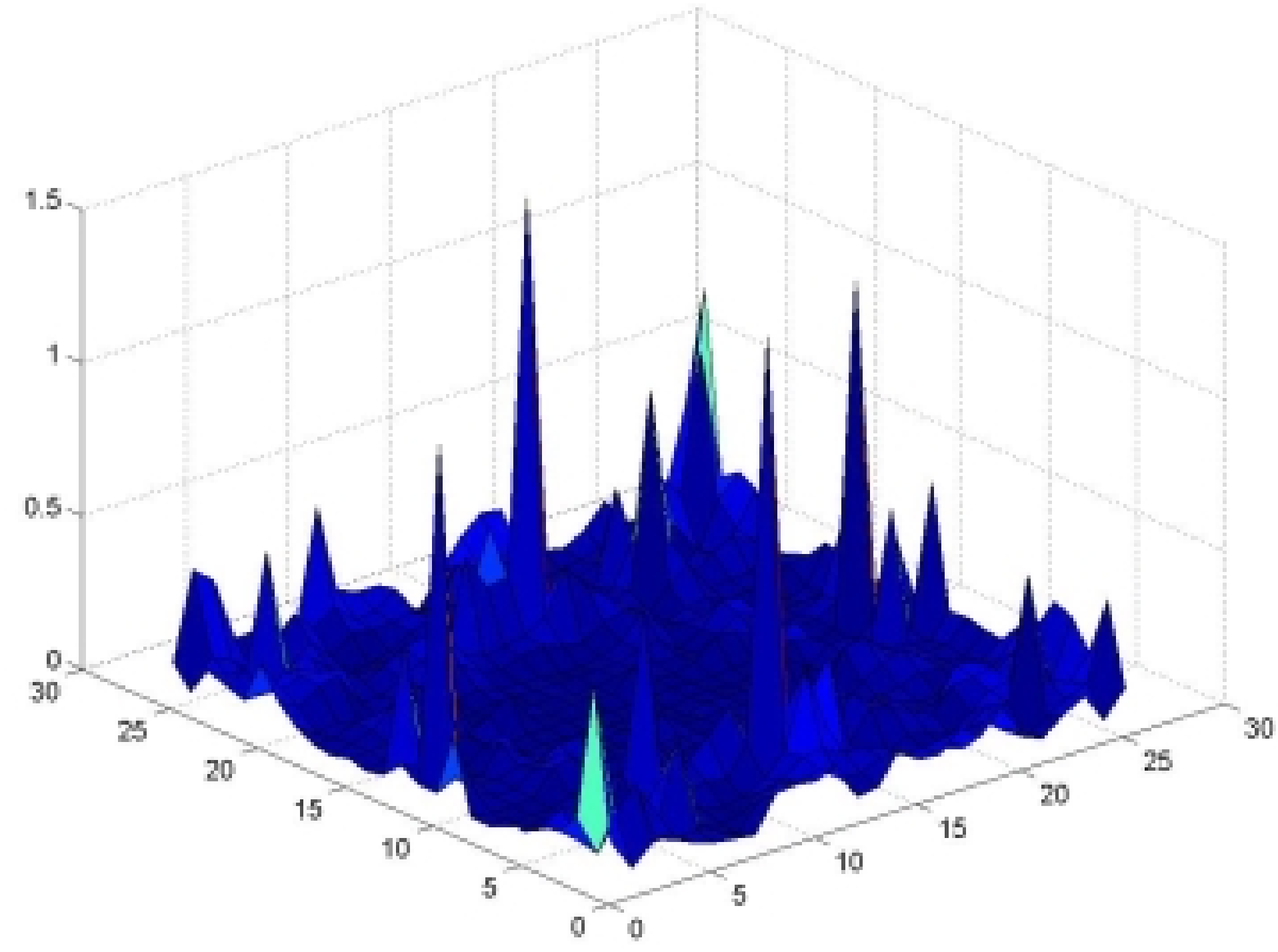

Figure 23

# REFERENCES

[1] J. P. Burg, *"Maximum Entropy Spectrum Analysis"*, **PH. D. DISSERTATION, STANFORD UNIV., 1975.**

[2] C. W. Therrien and H. T. El-Shaer, *"A Direct Algorithm for Computing 2-D AR Power Spectrum Estimates"*, **IEEE TRANS. ASSP, 37(1):106-117, JANUARY 1989.**

[3] C. W. Therrien, *"Relation Between 2-D and Multichannel Linear Prediction"*, IEEE TRANS. ASSP; VOL. ASSP-29, PP. 454-457, JUNE 1981.

[4] I. C. Gohberg and G. Heining, *"Inversion of finite Toeplitz matrices made of elements of a non-commutative algebra"*, REV. ROUMAINE MATH. PURES APPL., XIX(5):623-663, 1974, (IN RUSSIAN).

[5] I. C. Gohberg and A. Semencul, *"on the inversion of finite Toeplitz matrices and their continous analogs"*, MAT. ISSLED., 2(201:-233), 1972, (IN RUSSIAN).

[6] A. Jakobsson, S. L. Marple, P. Stoica, *"Computationally Efficient Two-Dimensional Capon Spectrum Analysis"*, IEEE TRANS. ON SIGNAL PROCESSING, VOL. 48, NO. 9, SEPTEMBER 2000.

[7] S.L Marple, Jr., *"Digital Spectral Analysis With Application"*, PRENTICE-HALL, ENGLEWOOD CLIFFS, N.J., 1987.

[8] G. Castro, *"Coefficients de Réflexion Généralisée. Extension de Covariance Multidimensionnelles et Autres Applications"*, PH.D. THESIS, UNIVERSITÉ DE PARIS-SUD CENTRE.

[9] Durbin J., "*The fitting of times series models*", REV. INST. STATIST., VOL.28 PP 233-244, (1960).

[10] Levinson N., *"The Wiener RMS (root-mean-square) error Criterion in filter design and prediction"*, J. MATH. ANAL. PHYS, VOL.25, PP 261-278, (1947).

[11] Wiener, N., "*Extrapolation, Interpolation, and smoothing of Stationary Time Series with Engenering Applications"*, JHON WILEY AND SONS, INC., NEW YORK, 1950.

[12] Kolmogorov, A., *"Sur l'interpolation et extrapolation des suites stationnaires"*, C.R. ACAD. SCI. PARIS, V. 208, 2043-2045, (1939).

[13] Shannon, Claude E. and Warren Weaver, *"The Mathematical Theory of Communication"*, THE UNIVERSITY IF ILLINOIS PRESS. URBANA. 1959.

[14] R. Kanhouche, "*Generalized Reflection Coefficients in Toeplitz-Block-Toeplitz Matrix Case and Fast Inverse 2D levinson Algorithm*", ***http://hal.ccsd.cnrs.fr/ccsd-00001420,*** APRIL 2004.

[15] Capon, J., *"High-Resolution-Frequency-Wavenumber Spectrum Analysis"*, **PROC. IEEE, VOL. 57, PP. 1408-1418, AUGUST 1969.**